\documentclass[leqno,10pt]{amsart}
 
\usepackage{amsmath,amsthm}
\usepackage{amssymb}
\usepackage{euscript}

\numberwithin{equation}{subsection}  

\newcommand{\sqsp}{\renewcommand{\baselinestretch}{1.1}\tiny\normalsize}

\newtheorem{thm}{Theorem}
\newtheorem{lemma}[thm]{Lemma}
\newtheorem{prop}[thm]{Proposition}
\newtheorem{cor}[thm]{Corollary}
\newtheorem{definition}[thm]{Definition}

\newcommand{\cat}[1]{{\EuScript #1}}
\newcommand{\cF}{\cat{F}}
\newcommand{\cC}{\cat{C}}

\newcommand{\del}{\partial}

\DeclareMathOperator{\Id}{Id} 
\DeclareMathOperator{\Hom}{\mathbf{Hom}}     
\DeclareMathOperator{\End}{\mathbf{End}}
\DeclareMathOperator{\Obs}{\mathbf{Obs}}
\DeclareMathOperator{\Endbar}{\overline{\mathbf{End}}}

\begin{document}
\title{Cohomology of $\lambda$-rings}
\author{Donald Yau}

\begin{abstract}
A cohomology theory for $\lambda$-rings is developed.  This is then applied to study deformations of $\lambda$-rings.
\end{abstract}

\email{dyau@math.uiuc.edu}
\address{Department of Mathematics, University of Illinois at Urbana-Champaign, 1409 W. Green Street, Urbana, IL 61801 USA}
\date{\today}

\maketitle

\sqsp

\section{Introduction}

The notion of a $\lambda$-ring was introduced by Grothendieck to study algebraic objects endowed with operations that act like exterior powers.  Since its introduction in the 1950s, $\lambda$-rings have been shown to play important roles in several areas of mathematics.  For example, in Algebraic Topology, the unitary $K$-theory of a topological space is a $\lambda$-ring.  When $X$ is a finite CW complex, the $\lambda$-operations on $K(X)$ are induced by exterior powers of vector bundles on $X$.  Similarly, the complex representation ring $R(G)$ of a group $G$ is a $\lambda$-ring with $\lambda$-operations given by the exterior powers of representations.  There is also an abundant supply of $\lambda$-rings from Algebra itself.  If $R$ is a commutative ring with unit, it can be shown that its universal Witt ring $\mathbf{W}(R)$ is always a $\lambda$-ring \cite{hazewinkel}.

The purposes of this note are (i) to introduce a cohomology theory for $\lambda$-rings and (ii) to use this to study $\lambda$-ring deformations along the lines of Gerstenhaber's theory \cite{ger}.

This note is organized as follows.  The following section contains a brief account of the basics of $\lambda$-rings and their Adams operations.  In Section \ref{sec:coh}, we define for a given $\lambda$-ring $R$ a cochain complex $\cF^*$ (see \eqref{eq:complex}) whose cohomology groups, denoted $H^*_\lambda(R)$, are the $\lambda$-ring cohomology groups of $R$.  Several basic observations are made.  First, the differential $d^n$ for $n \geq 1$ in $\cF^*$ is an alternating sum $\sum (-1)^i \del^i$.  There are ``codegeneracy'' maps $\sigma^i \colon \cF^n \to \cF^{n+1}$ for $n \geq 2$ such that the $\del^i$ and $\sigma^i$ satisfy the cosimplicial identities in dimensions $n \geq 2$.  In fact, $\cF^*$ is a subcomplex of a certain Hochschild cochain complex $\bar{\cF}^*$, defined in \S \ref{subsec:Fbar}, which coincides with $\cF^*$ in dimensions $2$ and above.  The cosimplicial identities in $\cF^*$ come from the cosimplicial abelian group that gives rise to the Hochschild complex $\bar{\cF}^*$.  Moreover, there is a composition product on $\cF^*$ that induces a product on cohomology, making $H^*_\lambda(R)$ a graded, associative, unital algebra (Corollary \ref{cor:product}).  The section ends with interpretations of $H^0_\lambda$ and $H^1_\lambda$ and the computation of these cohomology groups for the $\lambda$-ring $\mathbf{Z}$ (\S \ref{subsec:H}).

Section \ref{sec:def} is devoted to studying algebraic deformations of $\lambda$-rings, making use of the $\lambda$-ring cohomology in Section \ref{sec:coh}.  In particular, the infinitesimal deformation is a $1$-cocycle in $\cF^*$ whose cohomology class is well-defined by the equivalence class of the deformation (Proposition \ref{prop:inf}).  It follows that the vanishing of $H^1_\lambda(R)$ implies that $R$ is rigid (Corollary \ref{cor1:def}), meaning that every deformation of $R$ is equivalent to the trivial deformation.  The question of extending a $1$-cocycle to a deformation, or ``integrability'' in the terminology of Gerstenhaber \cite{ger}, is studied next.  Given a $1$-cocycle, the obstruction to extending it to a deformation is a sequence of $2$-cocycles (Theorem \ref{thm1:obs}).  This means that the simultaneous vanishing of their cohomology classes is equivalent to the extendibility of the given $1$-cocycle to a deformation (Corollary \ref{cor:int}).  It follows, in particular, that extendibility of a $1$-cocycle is automatic if $H^2_\lambda(R)$ is trivial (Corollary \ref{cor2:obs}).  The question of when two extensions are equivalent is also considered (Proposition \ref{prop:equivalence}).

One thing that is clearly missing in $\lambda$-ring cohomology is naturality.  A $\lambda$-ring map does not in general induce a map in $\lambda$-ring cohomology.  There is one exception, which is when the map is a $\lambda$-ring self-map.  This is due to the fact that the algebra of linear endomorphisms is used in the definition of $\lambda$-ring cohomology.  A map of rings, or even of $\lambda$-rings, does not in general induce a map on the algebras of linear endomorphisms.  So we only have naturality in the category whose sole object is the $\lambda$-ring under consideration and whose morphisms are its $\lambda$-ring self-maps.  Even in this restricted category, the induced map is only a map of graded groups, as it does not preserve the composition product.


\section{$\lambda$-rings and Adams operations}
\label{sec:lambda}

In preparation for studying $\lambda$-ring cohomology in the next two sections, in this section we briefly review some basic definitions about $\lambda$-rings and Adams operations.  For more discussions about $\lambda$-rings, consult Atiyah and Tall \cite{at} or Knutson \cite{knutson}.  The author's articles \cite{yau1,yau2} contain some recent results on $\lambda$-rings which might also be of interest to the reader.


\subsection{$\lambda$-rings}
\label{subsec:lambda}

By a $\lambda$-ring we mean a unital, commutative ring $R$ endowed with functions 
   \[
   \lambda^i \colon R ~\to~ R \quad (i \geq 0),
   \]
called $\lambda$-operations, which satisfy the following conditions.  For any integers $i, j \geq 0$ and elements $r$ and $s$ in $R$:
   \begin{itemize}
   \item $\lambda^0(r) = 1$.
   \item $\lambda^1(r) = r$.
   \item $\lambda^i(1) = 0$ for $i > 1$.
   \item $\lambda^i(r + s) = \sum_{k = 0}^i\, \lambda^k(r)\lambda^{i-k}(s)$.
   \item $\lambda^i(rs) = P_i(\lambda^1(r), \ldots , \lambda^i(r); \lambda^1(s), \ldots, \lambda^i(s))$.
   \item $\lambda^i(\lambda^j(r)) = P_{i,j}(\lambda^1(r), \ldots , \lambda^{ij}(r))$.
   \end{itemize}
The $P_i$ and $P_{i,j}$ are some universal polynomials with integer coefficients.  See the references mentioned above for the exact definitions of these polynomials.  Note that what we call a $\lambda$-ring here is sometimes called a ``special'' $\lambda$-ring in the literature.

For example, the ring of integers $\mathbf{Z}$ is a $\lambda$-ring with $\lambda^i(n) = \binom{n}{i}$.  In this case, all the Adams operations (to be reviewed below) are equal to the identity map on $\mathbf{Z}$.  This is the only $\lambda$-ring structure on $\mathbf{Z}$.

One important property of a $\lambda$-ring is that it must have characteristic $0$.  This can be seen from the linear map
   \[
   \lambda_t \colon R \to 1 + t R \lbrack \lbrack t \rbrack \rbrack 
   ~=~ \left\lbrace \sum \, a_i t^i \colon a_i \in R, \, a_0 = 1 \right\rbrace
   \]    
defined by
   \[
   \lambda_t(r) ~=~ \sum \, \lambda^i(r) t^i.
   \]
Here the additive group structure on $1 + t R \lbrack \lbrack t \rbrack \rbrack$ is given by the usual multiplication of power series.  The image of $n$ under $\lambda_t$ is $(1 + t)^n$, which is nonzero in $1 + t R \lbrack \lbrack t \rbrack \rbrack$ for any $n$.  In particular, for any positive integer $n$ and any prime $p$, the equation
   \[
   n^p ~\equiv~ n \pmod{pR}
   \]
holds.  This will be used in the next section when we study the $0$th $\lambda$-ring cohomology group.


\subsection{Adams operations}
\label{subsec:adams}

The $\lambda$-operations are sometimes hard to work with, since they are neither additive nor multiplicative.  One can extract ring maps from the  $\lambda$-operations, obtaining the so-called Adams operations
   \[
   \psi^n \colon R \to R \quad (n \geq 1).
   \]
More precisely, they are defined by the Newton formula:
   \[
   \psi^n(r) - \lambda^1(r)\psi^{n-1}(r) + \cdots + (-1)^{n-1}\lambda^{n-1}(r) \psi^1(r) + (-1)^n n\lambda^n(r) = 0.
   \]
The Adams operations satisfy the following properties:
\begin{itemize}
\label{adams}
\item All the $\psi^n$ are ring maps.
\item $\psi^1 = \Id$.
\item $\psi^m \psi^n = \psi^{mn} = \psi^n \psi^m$.
\item $\psi^p(r) \equiv r^p$ (mod $pR$) for each prime $p$ and element $r$ in $R$.
\end{itemize}

Suppose given a unital, commutative ring $R$ with self ring maps $\psi^n \colon R \to R$ satisfying the above four properties of Adams operations.  One can ask if it is possible to use the Newton formula to go backward and to produce a $\lambda$-ring structure on $R$.  This is, in fact, possible provided that $R$ is $\mathbf{Z}$-torsionfree.  More explicitly, a theorem of Wilkerson \cite{wil} says that if $R$ is as stated in the first sentence of this paragraph and is $\mathbf{Z}$-torsionfree, then there exists a unique $\lambda$-ring structure on $R$ whose Adams operations are exactly the given $\psi^n$.

We note that a ring $R$ with self ring maps $\psi^n \colon R \to R$ such that $\psi^1 = \Id$ and $\psi^m \psi^n = \psi^{mn}$ is sometimes called a ``weight system'' in the literature.  See, for example, Bar-Natan \cite{bn}.


\section{Cohomology of $\lambda$-rings}
\label{sec:coh}

The main purpose of this section is to introduce our $\lambda$-ring cohomology groups.  This is done in \ref{subsec:F}.  After that, we will discuss its connections with Hochschild cohomology in \ref{subsec:Fbar} and its product structure in \ref{subsec:product}.  The section closes with a discussion of the $0$th and the $1$st $\lambda$-ring cohomology groups.

Throughout this section, $R$ will denote a $\lambda$-ring with $\lambda$-operations $\lambda^i$ $(i \geq 0)$ and Adams operations $\psi^n$ $(n \geq 1)$.


\subsection{The complex $\cF^*$ and $\lambda$-ring cohomology}
\label{subsec:F}

To define the complex $\cF^* = \cF^*(R)$ that gives rise to $\lambda$-ring cohomology, we first need to establish some notations.

Denote by $\End(R)$ the (non-commutative) algebra of $\mathbf{Z}$-linear endomorphisms of $R$, in which the product is given by composition.  To make it clear that we are composing two endomorphisms $f$ and $g$, we will sometimes write $f \circ g$ instead of just $fg$.  We also need the following subalgebra of $\End(R)$.  Denote by $\Endbar(R)$ the subalgebra of $\End(R)$ consisting of those linear endomorphisms $f$ of $R$ that satisfy the condition,
   \begin{equation}
   \label{eq:endbar}
   f(r)^p ~\equiv~ f(r^p) \pmod{pR},
   \end{equation}
for every prime $p$ and each element $r \in R$.  We will use the symbol $T$ to denote the set of positive integers.

We are now ready to define the cochain complex $\cF^* = \cF^*(R)$.  Define $\cF^0$ to be the underlying additive group of $\Endbar(R)$ and $\cF^1$ to be the set of functions
   \[
   \label{eq:F1}
   f \colon T ~\to~ \End(R)
   \]
satisfying the condition, $f(p)(R) \subset pR$ for every prime $p$.  (The definitions of $\cF^0$ and $\cF^1$ might seem a little bit strange at first sight.  The reason for defining them as such will become apparent when we discuss deformations of $\lambda$-rings in the next section.)  For $n \geq 2$, $\cF^n$ is simply defined to be the set of functions
   \[
   f \colon T^n ~\to~ \End(R).
   \]
Each $\cF^n$ $(n \geq 1)$ inherits the obvious additive group structure from $\End(R)$.  Namely, if $f$ and $g$ are elements of $\cF^n$, then 
   \[
   (f + g)(m_1, \ldots, m_n)(r) ~=~ f(m_1, \ldots, m_n)(r) + g(m_1, \ldots, m_n)(r)
   \]
for $(m_1, \ldots, m_n) \in T^n$ and $r \in R$.

For $n \geq 0$, the differential $d^n \colon \cF^n \to \cF^{n+1}$ is defined by the formula
   \begin{multline}
   \label{eq:d}
   (d^n f)(m_0, \ldots, m_n) 
   ~=~ \psi^{m_0} \circ f(m_1, \ldots, m_n)
   ~+~ \sum_{i=1}^n \, (-1)^i f(m_0, \ldots, m_{i-1}m_i, \ldots, m_n) \\
   ~+~ (-1)^{n+1} f(m_0, \ldots, m_{n-1}) \circ \psi^{m_n}.
   \end{multline}
The $d^n$ are clearly additive group maps, and the only thing that we have to check is that the image of $d^0$ lies in $\cF^1$.  To see that this is the case, let $f$ be an element of $\cF^0$.  Then for any prime $p$ and element $r$ in $R$, we have that
   \[
   \begin{split}
   (d^0f)(p)(r) 
   & ~=~ \psi^p(f(r)) ~-~ f(\psi^p(r)) \\
   & ~\equiv~ f(r)^p ~-~ f(r^p) \pmod{pR} \\
   & ~\equiv~ 0 \pmod{pR}.
   \end{split}
   \]
This shows that $d^0$ is well-defined.

\medskip
\begin{lemma}
\label{lem:d}
For each $n \geq 0$, we have $d^{n+1} d^n = 0$.
\end{lemma}

\begin{proof}
The identity $d^1 d^0 = 0$ can be checked directly by writing out all six terms. 

For $n \geq 1$, we can write $d^n$ as $\sum_{i=0}^{n+1} (-1)^i \del^i$, where $\del^i \colon \cF^n \to \cF^{n+1}$ is the linear map given by
   \begin{equation}
   \label{eq:face map}
   (\del^i f)(m_0, \ldots, m_n) 
   ~=~ \begin{cases}
        \psi^{m_0} \circ f(m_1, \ldots, m_n) & \text{ if } i = 0 \\
        f(m_0, \ldots, m_{i-1}m_i, \ldots, m_n) & \text{ if } 1 \leq i \leq n \\
        f(m_0, \ldots, m_{n-1}) \circ \psi^{m_n} & \text{ if } i = n + 1. \end{cases}
   \end{equation}
Using the property, $\psi^n \psi^m = \psi^{mn}$, of the Adams operations, the ``cosimplicial identities" 
   \[
   \del^j \del^i ~=~ \del^i \del^{j-1} \quad (i < j)
   \]
can then be verified by direct inspection.  This implies, as usual, that $d^{n+1} d^n = 0$.
\end{proof}

Note that in this proof, we could have written $d^0$ formally as $\del^0 - \del^1$ just as above.  However, $\del^0$ and $\del^1$ do not necessarily have images in $\cF^1$.

The lemma gives us the cochain complex $\cF^* = \cF^*(R)$ of abelian groups,
   \begin{equation}
   \label{eq:complex}
   0 \to \cF^0 \xrightarrow{d^0} \cF^1 \xrightarrow{d^1} \cF^2 \xrightarrow{d^2} \cdots
   \end{equation}
with $\cF^n$ in dimension $n$.

\medskip
\begin{definition}
\label{def:coh}
The $n$th cohomology group of $\cF^* = \cF^*(R)$ is called the $n$th $\lambda$-\emph{ring cohomology group of} $R$, denoted by $H^n_\lambda(R)$.   
\end{definition}

The differentials $d^n$ look a lot like those in Hochschild cohomology theory.  There is, in fact, a close relationship between the complex $\cF^*$ and Hochschild theory, to which we now turn.


\subsection{Connections with Hochschild cohomology}
\label{subsec:Fbar}

Recall that $T$ denotes the set of positive integers and that $R$ is a $\lambda$-ring .  We will compare the complex $\cF^*$ with a certain Hochschild cochain complex.  For general discussions about Hochschild theory, refer, for example, to Weibel \cite{weibel}.

With the usual multiplication of integers, we can consider $T$ as a multiplicative, commutative monoid.  Then the underlying additive group of the algebra $\End(R)$ is a bimodule over the monoid-ring $\mathbf{Z}\lbrack T \rbrack$ via the action
   \[
   \begin{split}
   T \times \End(R) \times T & ~\to~ \End(R) \\
   (m,f,n)                   & ~\mapsto~ \psi^m \circ f \circ \psi^n,
   \end{split}
   \]
extended linearly to all of $\mathbf{Z} \lbrack T \rbrack$.  Therefore, we can consider the Hochschild cochain complex $C^* = C^*(\mathbf{Z}\lbrack T \rbrack, \End(R))$ of the monoid-ring $\mathbf{Z}\lbrack T \rbrack$ with coefficients in the bimodule $\End(R)$ and with ground ring $\mathbf{Z}$.  The $n$th cohomology group of $C^*$, denoted by $H^n(\mathbf{Z}\lbrack T \rbrack, \End(R))$, is called the $n$th Hochschild cohomology group of $\mathbf{Z}\lbrack T \rbrack$ with coefficients in $\End(R)$.

There is a canonical ring isomorphism
   \[
   \mathbf{Z}\lbrack T \rbrack^{\otimes n} 
   ~\cong~ \mathbf{Z} \lbrack T^n \rbrack,
   \]
where the multiplication on the monoid $T^n$ is defined coordinatewise.  Moreover, a $\mathbf{Z}$-linear map
   \[
   \mathbf{Z}\lbrack T^n \rbrack ~\to~ \End(R)
   \]
determines and is determined by a function 
   \[
   T^n ~\to~ \End(R).
   \]
Therefore, for $n \geq 2$, there is a canonical bijection
   \[
   C^n 
   ~=~ \Hom_\mathbf{Z}(\mathbf{Z}\lbrack T^n \rbrack, \End(R)) 
   ~\cong~ \cF^n
   \]
which, as one can check directly, respects the additive group structures.  Likewise, for $n = 0$ and $1$, one can identify $\cF^n$ canonically as a subgroup of $C^n$.  It is also straightforward to see from \eqref{eq:d} that, under the above identifications, the differentials in $\cF^*$ correspond to those in $C^*$.  This allows us to identify $\cF^*$ as a subcomplex of $C^*$, and the two complexes coincide from dimension $2$ onward.  In particular, we have the following result.

\medskip
\begin{prop}
\label{prop:Hoch}
There exist a canonical isomorphism
   \[
   H^n_\lambda(R) ~\cong H^n(\mathbf{Z}\lbrack T \rbrack, \End(R))
   \]
for each $n \geq 3$ and a canonical surjection
   \[
   H^2_\lambda(R) \twoheadrightarrow H^2(\mathbf{Z}\lbrack T \rbrack, \End(R)).
   \]
\end{prop}

It is well-known that the cochain complex $C^*$ arises from a cosimplicial abelian group $\cC^*$ with $C^n = \cC^n$ for all $n$.  In the proof of Lemma \ref{lem:d}, the maps $\del^i \colon \cF^n \to \cF^{n+1}$ for $n \geq 2$ are, under the identification $\cF^* \subset C^*$, exactly the coface maps of $\cC^*$.  There are also ``codegeneracy "maps 
   \[
   \sigma^i \colon \cF^{n+1} ~\to~ \cF^n \quad (i = 0, 1, \ldots, n)
   \]
defined by
   \[
   (\sigma^i f)(m_1, \ldots, m_n) ~=~ f(\ldots, m_i, 1, m_{i+1}, \ldots).
   \]
The maps $\del^i$ and $\sigma^i$ satisfy the usual cosimplicial identities in dimensions $2$ and above.  Once again, under the identification of $\cF^*$ as a subcomplex of $C^*$, these are the codegeneracy maps of $\cC^*$.  


\subsection{Composition product}
\label{subsec:product}

The purpose of this subsection is to observe that the $\lambda$-ring cohomology $H^*_\lambda(R)$ of $R$ is a graded ring.

\medskip
\begin{thm}
\label{thm:product}
Given a $\lambda$-ring $R$, there is an associative, bilinear pairing
   \[
   - \circ -  ~\colon~ \cF^n \otimes \cF^k ~\to~ \cF^{n+k} \quad (n, k \geq 0)
   \]
on the complex $\cF^* = \cF^*(R)$ with $\Id_R \in \cF^0$ as a two-sided identity.  This pairing satisfies the Leibnitz identity,
   \[
   d(f \circ g) ~=~ (df) \circ g ~+~ (-1)^{\vert f \vert}f \circ (dg),
   \]
where $\vert f \vert$ is the dimension of $f$.
\end{thm}

We call the pairing the \emph{composition product}.  The complex $\cF^*$ with the composition product is a differential graded algebra.  The Leibnitz identity implies that the product descends to cohomology with $\lbrack f \rbrack \circ \lbrack g \rbrack = \lbrack f \circ g \rbrack$, where $\lbrack f \rbrack$ denotes the cohomology class of a cocycle.

\medskip
\begin{cor}
\label{cor:product}
The composition product on $\cF^*$ induces a product on $H^*_\lambda(R)$, making it into a graded, associative, unital algebra.
\end{cor}

\begin{proof}[Proof of Theorem \ref{thm:product}]
The pairing is defined as follows.  Given $f \in \cF^n$, $g \in \cF^k$, and $(m_1, \ldots, m_{n+k}) \in T^{n+k}$, we set
   \[
   (f \circ g)(m_1, \ldots, m_{n+k}) 
   ~=~ f(m_1, \ldots, m_n) \circ g(m_{n+1}, \ldots, m_{n+k}),
   \]
where the $\circ$ on the right-hand side of the equation denotes composition of linear endomorphisms of $R$.  Associativity and bilinearity are straightforward to check, as is the assertion that $\Id_R$ acts as a two-sided identity.

As for the Leibnitz identity, let $f$ and $g$ be as above and let $(m_0, \ldots, m_{n+k})$ be in $T^{n+k+1}$.  Then $d(f \circ g)(m_0, \ldots, m_{n+k})$ is the sum of $n+k+2$ linear endomorphisms of $R$, $n+k$ of which come from the alternating sum $\sum (-1)^i (f \circ g)(\ldots, m_{i-1}m_i, \ldots)$.  Using the fact that
   \begin{multline*}
   (f \circ g)(\ldots, m_{i-1}m_i, \ldots) \\
   ~=~ \begin{cases}
       f(\ldots, m_{i-1}m_i, \ldots, m_n) \circ g(m_{n+1}, \ldots) & \text{ if } 1 \leq i \leq n \\
       f(m_0, \ldots, m_{n-1}) \circ g(\ldots, m_{i-1}m_i, \ldots) & \text{ if } n+1 \leq i \leq n+k,
       \end{cases}
   \end{multline*}
one observes that the terms for $1 \leq i \leq n$ (resp.\ $n + 1 \leq i \leq n + k$) correspond to the $n$ (resp.\ $k$) terms in $(df) \circ g$ (resp.\ $(-1)^{\vert f \vert} f \circ (dg)$) involving the alternating sum.  It follows easily from this observation that the Leibnitz identity holds.
\end{proof}

We remark that the composition product can also be defined on the Hochschild cochain complex $C^*$, and it has the same properties there.  Moreover, the subcomplex inclusion $\cF^* \subset C^*$ is a map of differential graded algebras, and the induced map on cohomology is a map of graded algebras.


\subsection{$H^0_\lambda$ and $H^1_\lambda$}
\label{subsec:H}

The purpose of this subsection is to discuss some basic properties of the $0$th and the $1$st $\lambda$-ring cohomology groups.  We will also compute these groups for the only $\lambda$-ring structure on $\mathbf{Z}$.

Recall that $\Endbar(R)$ is the group of linear endomorphisms $f$ of $R$ that satisfy the condition, $f(r)^p \equiv f(r^p) \pmod{pR}$, for each prime $p$ and each element $r \in R$.  The following result, which describes $H^0_\lambda$ explicitly, is immediate from the definition of $d^0$.

\medskip
\begin{prop}
\label{prop1:H0}
For any $\lambda$-ring $R$, we have that
   \[
   H^0_\lambda(R) 
   ~=~ \lbrace f \in \Endbar(R) ~\colon~ f \psi^n = \psi^n f \text{ for all n} \rbrace.
   \]
\end{prop}

Since a $\lambda$-ring $R$ must have characteristic $0$, for any integer $k$ and any element $r \in R$, the congruence relation
   \[
   (kr)^p ~\equiv k (r^p) \pmod{pR}
   \]
holds for each prime $p$.  This implies that the multiplication-by-$k$ endomorphism, $f_k \colon r \mapsto kr$, lies in $\Endbar(R)$.  It is also clear that this map commutes with $\psi^n$ for any $n$.  In particular, we have the following consequence of the proposition.

\medskip
\begin{cor}
\label{cor1:H0}
For any $\lambda$-ring $R$, $H^0_\lambda(R)$ contains $\mathbf{Z}$ as a canonical subgroup, which consists of the multiplication-by-$k$ endomorphisms of $R$.
\end{cor}

Recall that the ring of integers $\mathbf{Z}$ has a unique $\lambda$-ring structure given by $\lambda^i(n) = \binom{n}{i}$ with $\psi^m = \Id$ for all $m$.  Since any linear endomorphism $f$ of $\mathbf{Z}$ sends $n$ to $f(1)n$, a special case of the above corollary is

\medskip
\begin{cor}
\label{cor2:H0}
The $\lambda$-ring $\mathbf{Z}$ has $H^0_\lambda(\mathbf{Z}) \cong \mathbf{Z}$.
\end{cor}

We now turn to the group $H^1_\lambda$.

From the definition of $d^1$, the kernel of $d^1$ consists of those functions $f \in \cF^1$ such that
   \[
   f(mn) ~=~ \psi^m \circ f(n) ~~+ f(m) \circ \psi^n
   \]
for all $m$ and $n$.  Due to the similarity of this property with the defining property for derivations, we call these maps $\lambda$-\emph{derivations} (of $R$).  On the other hand, the image of $d^0$ consists of those functions $T \to \End(R)$ of the form
   \[
   \lbrack \psi^*, g \rbrack \colon n ~\longmapsto~ \psi^n \circ g - g \circ \psi^n
   \]
for some $g \in \Endbar(R)$.  In other words, they are just the functions obtained by ``twisting'' a $g \in \cF^0 = \Endbar(R)$ by $\psi^*$.  Because of this, we call these maps $\lambda$-\emph{inner derivations} (of $R$).

In particular, we have

\medskip
\begin{prop}
\label{prop:H1}
$H^1_\lambda(R)$ is the quotient of the group of $\lambda$-derivations by the group of $\lambda$-inner derivations.
\end{prop}

In the case of the $\lambda$-ring $\mathbf{Z}$, the Adams operations $\psi^n$ are all equal to the identity, so the only $\lambda$-inner derivation is $0$. On the other hand, identifying a linear endomorphism of $\mathbf{Z}$ with its image at $1$, one observes that a $\lambda$-derivation of $\mathbf{Z}$ is a function
   \[
   f \colon T ~\to~ \End(\mathbf{Z}) = \mathbf{Z}
   \]
such that $f(p) \in p\mathbf{Z}$ for each prime $p$ and that
   \[
   f(mn) ~=~ f(m) + f(n)
   \]
for all $m, n \geq 1$.  This second property simply means that, if $k$ has the prime factorization $p_1^{e_1} \cdots p_l^{e_l}$ with $e_i \geq 1$, then
   \[
   f(k) ~=~ e_1 f(p_1) + \cdots + e_l f(p_l).
  \]
In other words, the function $f$ is determined by the $f(p) \in p \mathbf{Z}$ for $p$ primes via this last equation.

Summarizing this discussion, we have

\medskip
\begin{cor}
\label{cor:H1}
$H^1_\lambda(\mathbf{Z}) \cong \prod_p\, p\mathbf{Z} \cong \prod_p \, \mathbf{Z}$, where the product is taken over the set of all primes.
\end{cor}


\section{Deformations of $\lambda$-rings}
\label{sec:def}

The purpose of this section is to study algebraic deformations of $\lambda$-rings along the path initiated by Gerstenhaber \cite{ger}, making use of the $\lambda$-ring cohomology developed in the previous section.

We remind the reader that $T$ denotes the set of positive integers and $R$ will always be an arbitrary $\lambda$-ring with Adams operations $\psi^n$.

Let us motivate the definition of a deformation of $R$ as follows.  Recall that the Adams operations are ring endomorphisms with the properties that $\psi^1 = \Id$, $\psi^{mn} = \psi^m \psi^n$, and $\psi^p(r) \equiv r^p \pmod{pR}$ for all primes $p$ and $r \in R$.  We would like to deform $R$ with respect to these properties.

Now let 
   \begin{equation}
   \label{eq:def}
   \Psi^*_t ~=~ \psi^*_0 + t\psi^*_1 + t^2\psi^*_2 + \cdots
   \end{equation}
be a formal power series, in which each $\psi^*_i$ is a function
   \[
   \psi^*_i \colon T ~\to~ \End(R)
   \]
with $\psi^*_0 = \psi^*$ (i.e.\ $\psi^n_0 = \psi^n$).  We will write $\psi^*_i(k)$ as $\psi^k_i$.  Then, in order for $\Psi^*_t$ to be a deformation of $R$, it should have the following properties:
\begin{itemize}
\item $\Psi^1_t = \Id$, meaning that
   \begin{equation}
   \label{eq:def1}
   \psi^1_i ~=~ 0 \quad (i \geq 1).
   \end{equation}

\item $\Psi^{mn}_t = \Psi^m_t \Psi^n_t$, meaning that
   \begin{equation}
   \label{eq:def2}
   \psi^{mn}_i ~=~ \sum_{j = 0}^i \, \psi^m_j \circ \psi^n_{i-j}
   \end{equation}
for all $i \geq 0$ and $m, n \geq 1$.

\item For each prime $p$, $\Psi^p_t(r) \equiv r^p \pmod{pR}$, which means that
   \begin{equation}
   \label{eq:def3}
   \psi^p_i(R) ~\subset~ pR \quad (i \geq 1).
   \end{equation}
In other words, $\psi^*_i \in \cF^1(R)$.
\end{itemize}

Observe that in \eqref{eq:def2}, if one takes $m = n = i = 1$, then the fact that $\psi^1 = \Id$ implies that $\psi^1_1 = 0$.  By an induction argument, still with $m = n = 1$, it follows that $\psi^1_i = 0$ for all $i \geq 1$.  In other words, \eqref{eq:def2} implies \eqref{eq:def1}, and we may disregard the latter.  We, therefore, define a \emph{deformation} of the $\lambda$-ring $R$ to be a formal power series $\Psi^*_t$ as in \eqref{eq:def} with each $\psi^*_i$ $(i \geq 1)$ in $\cF^1(R)$, satisfying the identity \eqref{eq:def2}.  Following Gerstenhaber \cite{ger}, the function $\psi^*_1$ is called the \emph{infinitesimal deformation} of $\Psi^*_t$.  In the rest of this section, we consider the following standard issues in algebraic deformation theory.

\begin{enumerate}
\item Identify the infinitesimal deformation with an appropriate cohomology class.
\item Obtain rigidity result from the previous step.
\item Describe cohomological obstructions to extending a cocycle to a deformation.
\item Describe cohomological obstructions to two such extensions being equivalent to each other.
\end{enumerate}

To do all this, we first need a suitable notion of equivalence of deformations.  Define a \emph{formal automorphism} of the $\lambda$-ring $R$ to be a formal power series
   \[
   \Phi_t ~=~ 1 + t \phi_1 + t^2 \phi_2 + \cdots,
   \]
in which each $\phi_i$ belongs to $\Endbar(R)$ with $1$ denoting the identity map on $R$.  Two deformations $\Psi^*_t$ and $\bar{\Psi}^*_t$ are said to be \emph{equivalent} if there exists a formal automorphism $\Phi_t$ such that
   \begin{equation}
   \label{eq:equiv}
   \bar{\Psi}^*_t ~=~ \Phi^{-1}_t \Psi^*_t \Phi_t.
   \end{equation}
This equation is to be understood in the following sense: if $f$ and $g$ are in $\Endbar(R)$, then $f \psi^*_i g$ is the function $T \to \End(R)$ given by
   \[
   (f \psi^*_i g)(n) ~=~ f \circ \psi^n_i \circ g.
   \]
It is straightforward to verify that if $\Psi^*_t$ is a deformation and if $\Phi_t$ is a formal automorphism, then $\bar{\Psi}^*_t$ defined by \eqref{eq:equiv} is also a deformation.

Now we can identify the infinitesimal deformation with a $1$-cocycle in $\cF^1(R)$.

\medskip
\begin{prop}
\label{prop:inf}
The infinitesimal deformation $\psi^*_1$ is a $1$-cocycle in the complex $\cF^*(R)$, and its cohomology class is well-defined by its equivalence class.
\end{prop}

\begin{proof}
The fact that $\psi^*_1$ is a $1$-cocycle follows directly from \eqref{eq:def2} (when $i = 1$).  If $\bar{\Psi}^*_t$ is a deformation that is equivalent to $\Psi^*_t$, then the difference $\bar{\psi}^*_1 - \psi^*_1$ is of the form $\lbrack \psi^*, \phi \rbrack$ for some $\phi \in \Endbar(R)$, and this is a $1$-coboundary.
\end{proof}

Suppose that in the deformation $\Psi^*_t$, one has $\psi^*_1 = \cdots = \psi^*_{l-1} = 0$ (i.e.\ $\psi^i_j = 0$ for all $i \geq 1$ and $j = 1, \ldots, l-1$).  Then one observes from \eqref{eq:def2} that $\psi^*_l$ is a $1$-cocycle.

\medskip
\begin{thm}
\label{thm1:def}
Suppose that $\Psi^*_t = \psi^* + t^l \psi^*_l + t^{l+1} \psi^*_{l+1} + \cdots$ is a deformation of a $\lambda$-ring $R$.  If $\psi^*_l$ is a $1$-coboundary in $\cF^1(R)$, then $\Psi^*_t$ is equivalent to a deformation of the form $\bar{\Psi}^*_t = \psi^* + t^{l+1}\bar{\psi}^*_{l+1} + t^{l+2}\bar{\psi}^*_{l+2} + \cdots$. 
\end{thm}

\begin{proof}
By assumption $\psi^*_l = \lbrack \psi^*, \phi_l \rbrack$ for some $\phi_l \in \Endbar(R)$.  Using the formal automorphism $\Phi_t = 1 - t^l \phi_l$, we see that $\Psi^*_t$ is equivalent to the deformation
   \[
   \begin{split}
   \bar{\Psi}^*_t
   &~=~ \Phi^{-1}_t \Psi^*_t \Phi_t \\
   &~\equiv~ (1 + t^l \phi_l)(\psi^* + t^l\psi^*_l)(1 - t^l \phi_l) \pmod{t^{l+1}} \\
   &~\equiv~ \psi^* + t^l(\psi^*_l - \lbrack \psi^*, \phi_l \rbrack) \pmod{t^{l+1}} \\
   &~\equiv~ \psi^* \pmod{t^{l+1}}.
   \end{split}
   \]
This finishes the proof.
\end{proof}

An immediate consequence of this result (and its proof) is a cohomological criterion for the rigidity of the $\lambda$-ring $R$.

\medskip
\begin{cor}
\label{cor1:def}
If $H^1_\lambda(R) = 0$, then every deformation of $R$ is equivalent to $\psi^*$.
\end{cor}

It was established in Proposition \ref{prop:inf} that the infinitesimal deformation is a $1$-cocycle in $\cF^*$.  This raises the question: Given a $1$-cocycle, is it the infinitesimal deformation of a deformation?  To what extent is this deformation unique?   We will break each one of these questions into a sequence of ``smaller" questions, which we then approach from an obstruction-theoretic view point.

Fix a $\lambda$-ring $R$.  Following Gerstenhaber and Wilkerson \cite{gw}, we define, for each $N \geq 1$, a \emph{deformation of order} $N$ to be a formal power series
   \[
   \Psi^*_t ~=~ \psi^* + t\psi^*_1 + \cdots + t^N \psi^*_N
   \]
with each $\psi^*_i$ in $\cF^1(R)$, satisfying the identity \eqref{eq:def2} modulo $t^{N+1}$.  This last requirement simply means that 
   \[
   \Psi^m_t \Psi^n_t ~=~ \Psi^{mn}_t \pmod{t^{N+1}}
   \]
for all $m, n \geq 1$.  One can think of a deformation as a deformation of order $\infty$.  A formal automorphism is defined just as before, and two deformations of order $N$ are said to be equivalent if there exists a formal automorphism for which \eqref{eq:equiv} holds modulo $t^{N+1}$.  We say that $\Psi^*_t$ \emph{extends to order} $N + 1$ if there exists an element $\psi^*_{N+1} \in \cF^1(R)$ such that the formal power series
   \begin{equation}
   \label{eq:extension}
   \bar{\Psi}^*_t ~=~ \Psi^*_t + t^{N+1} \psi^*_{N+1}
   \end{equation}
is a deformation of order $N + 1$.   We call $\bar{\Psi}^*_t$ an \emph{order} $N + 1$ \emph{extension of} $\Psi^*_t$.

Let $\Psi^*_t$ be a deformation of order $N$.  Consider the function 
   \[
   \Obs(\Psi^*_t) \colon T^2 ~\to~ \End(R)
   \]
defined by
   \[
   \Obs(\Psi^*_t)(m,n) ~=~ - \sum_{i=1}^N\, \psi_i^m \circ \psi^n_{N+1-i}.
   \]

\medskip
\begin{lemma}
\label{lem1:obs}
The element $\Obs(\Psi^*_t) \in \cF^2(R)$ is a $2$-cocycle.
\end{lemma}

\begin{proof}
Recall that we can write $d^2 \colon \cF^2 \to \cF^3$ as $\sum_{i=0}^3\, (-1)^i \del^i$ (see \eqref{eq:face map}).  For any triple $(m_0, m_1, m_2) \in T^3$, we have  
   \[
   \begin{split}
   (\del^0 - \del^1)(\Obs(\Psi^*_t))(m_0,m_1,m_2) 
   &~=~ \sum_{\substack{i+j+k = N+1 \\ i,k>0}} \, \psi^{m_0}_i \circ \psi^{m_1}_j \circ \psi^{m_2}_k \\
   &~=~ -(\del^2 - \del^3)((\Obs(\Psi^*_t))(m_0,m_1,m_2).
   \end{split}
   \]
It follows that $d^2 \Obs(\Psi^*_t) = 0$.
\end{proof}

Now suppose that $\psi^*_{N+1}$ is an element of $\cF^1(R)$.  Consider the formal power series $\bar{\Psi}^*_t$ defined by \eqref{eq:extension}.  It is an order $N + 1$ extension of $\Psi^*_t$ if and only if \eqref{eq:def2} holds when $i = N + 1$.  This is true, since the identities for $i \leq N$ in \eqref{eq:def2} automatically hold, as they only involve $\psi^*_i$ for $i \leq N$.  Collecting the three terms in \eqref{eq:def2} (with $i = N+1$) involving $\psi^*_{N+1}$, \eqref{eq:def2} can be rewritten as
   \begin{equation}
   \label{eq1:obs}
   (d^1 \psi^*_{N+1})(m,n) ~=~ \Obs(\Psi^*_t)(m,n).
   \end{equation}
In other words, $\Psi^*_t$ extends to order $N+1$ if and only if $\Obs(\Psi^*_t)$ is a $2$-coboundary.

Summarizing, we have determined the obstructions to extending a deformation of order $N$ to a deformation of one higher order.  We record it as follows.

\medskip
\begin{thm}
\label{thm1:obs}
Let $\Psi^*_t$ be a deformation of order $N$.  Then it extends to order $N + 1$ if and only if the $2$-cocycle $\Obs(\Psi^*_t)$ is cohomologous to $0$.
\end{thm}

Starting with a $1$-cocycle, we obtain the obstructions to extending it to a deformation by applying this theorem repeatedly.

\medskip
\begin{cor}
\label{cor:int}
Let $\psi^*_1 \in \cF^1(R)$ be a $1$-cocycle.  Then there exists a sequence of (obstruction) classes $\omega_i$ $(i = 1, 2, \ldots)$ in $H^2_\lambda(R)$, where $\omega_n$ is defined if and only if $\omega_1, \ldots, \omega_{n-1}$ are all defined and equal to $0$.  Moreover, the deformation of order $1$, $\Psi^*_t = \psi^* + t\psi^*_1$, extends to a deformation if and only if $\omega_i$ is defined and equal to $0$ for each $i = 1, 2, \ldots$. 
\end{cor}

In particular, we have the following cohomological condition that guarantees the existence of extensions.

\medskip
\begin{cor}
\label{cor2:obs}
If $H^2_\lambda(R) = 0$, then every deformation of order $N \geq 1$ extends to a deformation.
\end{cor}

Finally, we consider the question of whether two extensions are equivalent.  Let $\Psi^*_t$ be a deformation of order $N$ and let $\bar{\Psi}^*_t$ and $\tilde{\Psi}^*_t$ be two order $N + 1$ extensions of $\Psi^*_t$.  Then it follows from the way the obstruction class is defined that
   \[
   d^1 \bar{\psi}^*_{N+1} ~=~ \Obs(\Psi^*_t) ~=~ d^1 \tilde{\psi}^*_{N+1}.
   \]
In particular, $\bar{\psi}^*_{N+1} - \tilde{\psi}^*_{N+1}$ is a $1$-cocycle.

\medskip
\begin{prop}
\label{prop:equivalence}
If the $1$-cocycle $\bar{\psi}^*_{N+1} - \tilde{\psi}^*_{N+1}$ is cohomologous to $0$, then the two deformations $\bar{\Psi}^*_t$ and $\tilde{\Psi}^*_t$ of order $N + 1$ are equivalent.
\end{prop}

\begin{proof}
By assumption we have $\bar{\psi}^*_{N+1} - \tilde{\psi}^*_{N+1} = \lbrack \psi^*, \phi \rbrack$ for some $\phi \in \Endbar(R)$.  Define the formal automorphism $\Phi_t = 1 + t^{N+1}\phi$.  Then an argument essentially identical to the proof of Theorem \ref{thm1:def} shows that $\bar{\Psi}^*_t \equiv \Phi^{-1}_t \tilde{\Psi}^*_t \Phi_t \pmod{t^{N+2}}$.
\end{proof}

The author is not sure whether the converse of this proposition is true or not.

\section*{Acknowledgment}
The author would like to thank the referee for his/her careful reading of the manuscript.  He would also like to thank James McClure for introducing the subject of algebraic deformations to him, and Charles Rezk for discussions related to this project.



\end{document}